# Analyzing and enhancing OSKI for sparse matrix-vector multiplication [1]

Kadir Akbudak[a], Enver Kayaaslan[a], Cevdet Aykanat[a,*]

[a]Computer Engineering Department, Bilkent University, Ankara, Turkey

**Abstract**

Sparse matrix-vector multiplication (SpMxV) is a kernel operation widely used in iterative linear solvers. The same sparse matrix is multiplied by a dense vector repeatedly in these solvers. Matrices with irregular sparsity patterns make it difficult to utilize cache locality effectively in SpMxV computations. In this work, we investigate single- and multiple-SpMxV frameworks for exploiting cache locality in SpMxV computations. For the single-SpMxV framework, we propose two cache-size-aware top-down row/column-reordering methods based on 1D and 2D sparse matrix partitioning by utilizing the column-net and enhancing the row-column-net hypergraph models of sparse matrices. The multiple-SpMxV framework depends on splitting a given matrix into a sum of multiple nonzero-disjoint matrices so that the SpMxV operation is performed as a sequence of multiple input- and output-dependent SpMxV operations. For an effective matrix splitting required in this framework, we propose a cache-size-aware top-down approach based on 2D sparse matrix partitioning by utilizing the row-column-net hypergraph model. The primary objective in all of the three methods is to maximize the exploitation of temporal locality. We evaluate the validity of our models and methods on a wide range of sparse matrices by performing actual runs through using OSKI. Experimental results show that proposed methods and models outperform state-of-the-art schemes.

## 1. Introduction

Sparse matrix-vector multiplication (SpMxV) is an important kernel operation in iterative linear solvers used for the solution of large, sparse, linear systems of equations. In these iterative solvers, the SpMxV operation $y \leftarrow Ax$ is repeatedly performed with the same large, irregularly sparse matrix $A$. Irregular access pattern during these repeated SpMxV operations causes poor usage of CPU caches in today's deep memory hierarchy technology. However, SpMxV operation has a potential to exhibit very high performance gains if temporal and spatial localities are respected and exploited properly. Here, temporal locality refers to the reuse of data words (e.g., $x$-*vector* entries) within relatively small time durations, whereas spatial locality refers to the use of data words (e.g., matrix nonzeros) within relatively close storage locations (e.g., in the same lines).

In this work, we investigate two distinct frameworks for the SpMxV operation: single-SpMxV and multiple-SpMxV frameworks. In the single-SpMxV framework, the $y$-*vector* results are computed by performing a single SpMxV operation $y \leftarrow Ax$. In the multiple-SpMxV framework, $y \leftarrow Ax$ operation is computed as a sequence of multiple input- and output-dependent SpMxV operations, $y \leftarrow y + A^k x$ for $k = 1, \ldots, K$, where $A = A^1 + \cdots + A^K$. For the single-SpMxV framework, we propose two cache-size-aware row/column reordering methods based on top-down 1D and 2D partitioning of a given sparse matrix. The 1D-partitioning-based method relies on transforming a sparse matrix into a singly-bordered block-diagonal (SB) form by utilizing the column-net hypergraph model [4–6]. The 2D-partitioning-based method relies on transforming a sparse matrix into a doubly-bordered block-diagonal (DB) form by utilizing the row-column-net hypergraph model [4,8]. We provide upper bounds on the number of cache misses based on these transformations, and show that the objectives in the transformations based on partitioning the respective hypergraph models correspond to minimizing these upper bounds. In the 1D-partitioning-based method, the column-net hypergraph model correctly encapsulates the minimization of the respective upper bound. For the 2D-partitioning-based method, we propose an enhancement to the row-column-net hypergraph model to encapsulate the minimization of the respective upper bound on the number of cache misses. The primary objective in both methods is to maximize the exploitation of the temporal locality due to the access of $x$-*vector* entries, whereas exploitation of the spatial locality due to the access of $x$-$vector$ entries is a secondary objective. In this paper, we claim that exploiting temporal locality is more important than exploiting spatial locality (for practical line sizes) in SpMxV operations that involve irregularly sparse matrices.

---

[*]Corresponding author.
   tel. +90-312-290-1625   fax. +90-312-266-4047   e-mail. aykanat@cs.bilkent.edu.tr
[1]A complete version of this work was submitted to *SIAM Journal on Scientific Computing (SISC)*.

The multiple-SpMxV framework depends on splitting a given matrix into a sum of multiple nonzero-disjoint matrices so that the SpMxV operation is computed as a sequence of multiple SpMxV operations. For an effective matrix splitting required in this framework, we propose a cache-size-aware top-down approach based on 2D sparse matrix partitioning by utilizing the row-column-net hypergraph model [4, 8]. We provide an upper bound on the number of cache misses based on this matrix-splitting, and show that the objective in the hypergraph-partitioning (HP) based matrix partitioning exactly corresponds to minimizing this upper bound. The primary objective in this method is to maximize the exploitation of the temporal locality due to the access of both $x$-$vector$ and $y$-$vector$ entries.

We evaluate the validity of our models and methods on a wide range of sparse matrices. The experiments are carried out by performing actual runs through using OSKI (BeBOP Optimized Sparse Kernel Interface Library) [17]. Experimental results show that the proposed methods and models outperform state-of-the-art schemes and also these results conform to our expectation that temporal locality is more important than spatial locality in SpMxV operations that involve irregularly sparse matrices.

The rest of the paper is organized as follows: Background material is introduced in Section 2. The two frameworks along with our contributed models and methods are described in Sections 3 and 4. We present the experimental results in Section 5. Finally, the paper is concluded in Section 6.

## 2. Background

### 2.1. Sparse-matrix storage schemes

There are two standard sparse-matrix storage schemes for the SpMxV operation: *Compressed Storage by Rows* (CSR) and *Compressed Storage by Columns* (CSC) [10, 13]. In this paper, we restrict our focus on the SpMxV operation using the CSR storage scheme without loss of generality, whereas cache aware techniques such as prefetching, blocking, etc. are out of the scope of this paper. In the following paragraphs, we review the standard CSR scheme and a CSR variant.

The CSR scheme contains three 1D arrays: *nonzero*, *colIndex* and *rowStart*. The values and the column indices of nonzeros are respectively stored in row-major order in the *nonzero* and *colIndex* arrays in a one-to-one manner. The *rowStart* array stores the index of the first nonzero of each row in the *nonzero* and *colIndex* arrays.

The *Incremental Compressed Storage by Rows* (ICSR) scheme [12] is reported to decrease instruction overhead by using pointer arithmetic. In ICSR, the *colIndex* array is replaced with the *colDiff* array, which stores the increments in the column indices of the successive nonzeros stored in the nonzero array. The *rowStart* array is replaced with the *rowJump* array which stores the increments in the row indices of the successive nonzero rows. The ICSR scheme has the advantage of handling zero rows efficiently since it avoids the use of the *rowStart* array. This feature of ICSR is exploited in our multiple-SpMxV framework since this scheme introduces many zero rows in the individual sparse matrices. Details of the SpMxV algorithms utilizing CSR and ICSR are described in our technical report [1].

### 2.2. Data locality in CSR-based SpMxV

In accessing matrix nonzeros, temporal locality is not feasible since the elements of each of the *nonzero*, *colIndex* (*colDiff* in ICSR) and *rowStart* (*rowJump* in ICSR) arrays are accessed only once. Spatial locality is feasible and it is achieved automatically by nature of the CSR scheme since the elements of each of these three arrays are accessed consecutively.

In accessing $y$-$vector$ entries, temporal locality is not feasible since each $y$-$vector$ result is written only once to the memory. As a different view, temporal locality can be considered as feasible but automatically achieved especially at the register level because of the summation of scalar nonzero and $x$-$vector$ entry product results to the temporary variable. Spatial locality is feasible and it is achieved automatically since the $y$-$vector$ entry results are stored consecutively.

In accessing $x$-$vector$ entries, both temporal and spatial localities are feasible. Temporal locality is feasible since each $x$-$vector$ entry may be accessed multiple times. However, exploiting the temporal and spatial localities for the $x$-$vector$ is the major concern in the CSR scheme since $x$-$vector$ entries are accessed through a *colIndex* array (*colDiff* in ICSR) in a non-contiguous and irregular manner.

## 3. Single-SpMxV framework

In this framework, the $y$-$vector$ results are computed by performing a single SpMxV operation, i.e., $y \leftarrow Ax$. The objective in this scheme is to reorder the columns and rows of matrix $A$ for maximizing the exploitation of temporal and spatial locality in accessing $x$-$vector$ entries. That is, the objective is to find row and column permutation matrices $P_r$ and $P_c$ so that $y \leftarrow Ax$ is computed as $\hat{y} \leftarrow \hat{A}\hat{x}$, where $\hat{A} = P_r A P_c$, $\hat{x} = x P_c$ and $\hat{y} = P_r y$. For the sake of simplicity of presentation, reordered input and output vectors $\hat{x}$ and $\hat{y}$ will be referred to as $x$ and $y$ in the rest of the paper.

Recall that temporal locality in accessing $y$-*vector* entries is not feasible, whereas spatial locality is achieved automatically because $y$-*vector* results are stored and processed consecutively. Reordering the rows with similar sparsity pattern nearby increases the possibility of exploiting temporal locality in accessing $x$-*vector* entries. Reordering the columns with similar sparsity pattern nearby increases the possibility of exploiting spatial locality in accessing $x$-*vector* entries. This row/column reordering problem can also be considered as a row/column clustering problem and this clustering process can be achieved in two distinct ways: top-down and bottom-up. In this section, we propose and discuss cache-size-aware top-down approaches based on 1D and 2D partitioning of a given matrix. Although a bottom-up approach based on hierarchical clustering of rows/columns with similar patterns is feasible, such a scheme is not discussed in this work.

In Sections 3.1 and 3.2, we present two theorems that give the guidelines for a "good" cache-size-aware row/column reordering based on 1D and 2D matrix partitioning. These theorems provide upper bounds on the number of cache misses due to the access of $x$-*vector* entries in the SpMxV operation performed on sparse matrices in two special forms, namely SB and DB forms. In these theorems, $\Phi_x(A)$ denotes the number of cache misses due to the access of $x$-*vector* entries in a CSR-based SpMxV operation to be performed on matrix $A$.

In the theorems given in Sections 3 and 4, fully associative cache is assumed, since misses in a fully associative cache are capacity misses and are not conflict misses. In these theorems, a matrix/submatrix is said to fit into the cache if the size of the CSR storage of the matrix/submatrix together with the associated $x$ and $y$ vectors/subvectors is smaller than the size of the cache. The proofs of these theorems are provided in our technical report [1].

*3.1. Row/column reordering based on 1D matrix partitioning*

We consider a row/column reordering which permutes a given matrix $A$ into a $K$-way columnwise SB form

$$\hat{A} = A_{SB} = P_r A P_c = \begin{bmatrix} A_{11} & & & & A_{1B} \\ & A_{22} & & & A_{2B} \\ & & \ddots & & \vdots \\ & & & A_{KK} & A_{KB} \end{bmatrix} = \begin{bmatrix} R_1 \\ R_2 \\ \vdots \\ R_K \end{bmatrix}$$

$$= \begin{bmatrix} C_1 & C_2 & \ldots & C_K & C_B \end{bmatrix}. \qquad (1)$$

Here, $A_{kk}$ denotes the $k$th diagonal block of $A_{SB}$. $R_k = [0 \ldots 0 \; A_{kk} \; 0 \ldots 0 \; A_{kB}]$ denotes the $k$th row slice of $A_{SB}$, for $k = 1, \ldots, K$. $C_k = \begin{bmatrix} 0 \ldots 0 \; A_{kk}^T \; 0 \ldots 0 \end{bmatrix}^T$ denotes the $k$th column slice of $A_{SB}$, for $k = 1, \ldots, K$, and $C_B$ denotes the column border as follows

$$C_B = \begin{bmatrix} A_{1B} \\ A_{2B} \\ \vdots \\ A_{KB} \end{bmatrix}. \qquad (2)$$

Each column in the border $C_B$ is called a *row-coupling column* or simply a *coupling column*. Let $\lambda(c_j)$ denote the number of $R_k$ submatrices that contain at least one nonzero of column $c_j$ of matrix $A_{SB}$, i.e.,

$$\lambda(c_j) = |\{R_k \in A_{SB} : c_j \in R_k\}|. \qquad (3)$$

In this notation, a column $c_j$ is a coupling column if $\lambda(c_j) > 1$. Here and hereafter, a submatrix notation is interchangeably used to denote both a submatrix and the set of non-empty rows/columns that belong to that matrix. For example, in (3), $R_k$ denotes both the $k$th row slice of $A_{SB}$ and the set of columns that belong to submatrix $R_k$.

The individual $y \leftarrow Ax$ can be equivalently represented as $K$ output-independent but input-dependent SpMxV operations, i.e., $y_k \leftarrow R_k x$ for $k = 1, \ldots, K$, where each submatrix $R_k$ is assumed to be stored in CSR scheme. These SpMxV operations are input dependent because of the $x$-*vector* entries corresponding to the coupling columns.

**Theorem 1** *Given a $K$-way SB form $A_{SB}$ of matrix $A$ such that every submatrix $R_k$ fits into the cache, then we have*

$$\Phi_x(A_{SB}) \leq \sum_{c_j \in A_{SB}} \lambda(c_j) \qquad (4)$$

Theorem 1 leads us to a cache-size-aware top-down row/column reordering through an $A$-to-$A_{SB}$ transformation that minimizes the upper bound given in (4) for $\Phi_x(A_{SB})$. Minimizing this sum relates to minimizing the number of cache misses due to the loss of temporal locality.

As discussed in [2], this $A$-to-$A_{SB}$ transformation problem can be formulated as an HP problem using the column-net model of matrix $A$ with the part size constraint of cache size and the partitioning objective of minimizing cutsize according to the connectivity metric definition given in [1]. In this way, minimizing the cutsize corresponds to minimizing the upper bound given in Theorem 1 for the number of cache misses due to the access of $x$-*vector* entries. This reordering method will be referred to as "sHP$_{\text{CN}}$", where the small letter "s" is used to indicate the single-SpMxV framework.

*3.2. Row/column reordering based on 2D matrix partitioning*

We consider a row/column reordering which permutes a given matrix $A$ into a $K$-way DB form

$$\hat{A} = A_{DB} = P_r A P_c = \begin{bmatrix} A_{11} & & & & A_{1B} \\ & A_{22} & & & A_{2B} \\ & & \ddots & & \vdots \\ & & & A_{KK} & A_{KB} \\ A_{B1} & A_{B2} & \ldots & A_{BK} & A_{BB} \end{bmatrix} = \begin{bmatrix} R_1 \\ R_2 \\ \vdots \\ R_K \\ R_B \end{bmatrix} = \begin{bmatrix} A'_{SB} \\ R_B \end{bmatrix}$$

$$= \begin{bmatrix} C_1 & C_2 & \ldots & C_K & C_B \end{bmatrix}. \tag{5}$$

Here, $R_B = [A_{B1} \ A_{B2} \ \ldots \ A_{BK} \ A_{BB}]$ denotes the row border. Each row in $R_B$ is called a *column-coupling row* or simply a *coupling row*. $A'_{SB}$ denotes the columnwise SB part of $A_{DB}$ excluding the row border $R_B$. $R_k$ denotes the $k$th row slice of both $A'_{SB}$ and $A_{DB}$. $\lambda'(c_j)$ denotes the connectivity of column $c_j$ in $A'_{SB}$. $C'_B$ denotes the column border of $A'_{SB}$, whereas $C_B = [C'^T_B \ A^T_{BB}]^T$ denotes the column border of $A_{DB}$. $C_k = \begin{bmatrix} 0 \ldots 0 \ A^T_{kk} \ 0 \ldots 0 \ A^T_{Bk} \end{bmatrix}^T$ denotes the $k$th column slice of $A_{DB}$.

**Theorem 2** *Given a K-way DB form $A_{DB}$ of matrix $A$ such that every submatrix $R_k$ of $A'_{SB}$ fits into the cache, then we have*

$$\Phi_x(A_{DB}) \leq \sum_{c_j \in A'_{SB}} \lambda'(c_j) + \sum_{r_i \in R_B} nnz(r_i). \tag{6}$$

Theorem 2 leads us to a cache-size-aware top-down row/column reordering through an $A$-to-$A_{DB}$ transformation that minimizes the upper bound given in (6) for $\Phi_x(A_{DB})$. Here, minimizing this sum relates to minimizing the number of cache misses due to the loss of temporal locality.

Here we propose to formulate the above-mentioned $A$-to-$A_{DB}$ transformation problem as an HP problem using the row-column-net model of matrix $A$ with a part size constraint of cache size. In the proposed formulation, column nets are associated with unit cost (i.e., $cost(n^c_j) = 1$ for each column $c_j$ and the cost of each row net is set to the number of nonzeros in the respective row (i.e., $cost(n^r_i) = nnz(r_i)$). However, existing HP tools do not handle a cutsize definition that encapsulates the right-hand side of (6), because the connectivity metric should be enforced for column nets, whereas the cut-net metric should be enforced for row nets. In order to encapsulate this different cutsize definition, we adapt and enhance the cut-net removal and cut-net splitting techniques adopted in RB algorithms utilized in HP tools. The details of the enhanced row-column-net model can be found in our technical report [1].

## 4. Multiple-SpMxV framework

Let $\Pi = \{A^1, A^2, \ldots, A^K\}$ denote a splitting of matrix $A$ into $K$ $A^k$ matrices, where $A = A^1 + A^2 + \cdots + A^K$. In $\Pi$, $A^k$ matrices are mutually nonzero-disjoint, however they are not necessarily row disjoint or column disjoint. Note that every splitting $\Pi$ defines an access order on the matrix non-zeros, and every access order can define $\Pi$ that causes it.

In this framework, $y \leftarrow Ax$ operation is computed as a sequence of $K$ input- and output-dependent SpMxV operations, $y \leftarrow y + A^k x$ for $k = 1, \ldots, K$. Individual SpMxV results are accumulated in the output vector $y$ on the fly in order to avoid additional write operations. The individual SpMxV operations are input dependent because of the shared columns among the $A^k$ matrices, whereas they are output dependent because of the shared rows among the $A^k$ matrices. Note that $A^k$ matrices are likely to contain empty rows and columns. The splitting of matrix $A$ should be done in such a way that the temporal and spatial localities of individual SpMxVs are exploited in order to minimize the number of cache misses.

In Section 4.1, we present a theorem that gives the guidelines for a "good" cache-size-aware matrix splitting based on 2D matrix partitioning. This theorem provides an upper bound on the total number of cache misses due to the access of $x$-vector and $y$-vector entries in all $y \leftarrow y + A^k x$ operations.

*4.1. Splitting $A$ into $A^k$ matrices based on 2D matrix partitioning*

Given a splitting $\Pi$ of matrix $A$, let $\Phi_x(A, \Pi)$ and $\Phi_y(A, \Pi)$ respectively denote the number of cache misses due to the access of $x$-vector and $y$-vector entries during $y \leftarrow y + A^k x$ operations for $k = 1, \ldots, K$. Here, the total number of cache misses can be expressed as $\Phi(A, \Pi) = \Phi_x(A, \Pi) + \Phi_y(A, \Pi)$. Let $\lambda(r_i)$ and $\lambda(c_j)$ respectively denote the number of $A^k$ matrices that contain at least one nonzero of row $r_i$ and one nonzero of column $c_j$ of matrix $A$, i.e., $\lambda(r_i) = |\{A^k \in \Pi : r_i \in A^k\}|$ and $\lambda(c_j) = |\{A^k \in \Pi : c_j \in A^k\}|$.

**Theorem 3** *Given a splitting $\Pi = \{A^1, A^2, \ldots, A^K\}$ of matrix $A$, then we have*
  a) $\Phi_x(A, \Pi) \leq \sum_{c_j \in A} \lambda(c_j)$, *if each $A^k$ matrix fits into the cache;*
  b) $\Phi_y(A, \Pi) \leq \sum_{r_i \in A} \lambda(r_i)$.

**Corollary 4** *If each $A^k$ in $\Pi$ fits into the cache, then we have*

$$\Phi(A, \Pi) \leq \sum_{r_i \in A} \lambda(r_i) + \sum_{c_j \in A} \lambda(c_j). \tag{7}$$

Corollary 4 leads us to a cache-size-aware top-down matrix splitting that minimizes the upper bound given in (7) for $\Phi(A,\Pi)$. Here, minimizing this sum relates to minimizing the number of cache misses due to the loss of temporal locality.

The matrix splitting problem can be formulated as an HP-based 2D matrix paritifioning using the row-column-net model [4, 8] of matrix $A$ with a part size constraint of cache size and partitioning objective of minimizing cutsize according to the connectivity metric definition given in [1]. In this way, minimizing the cutsize corresponds to minimizing the upper bound given in Theorem 3 for the total number of cache misses due to the access of *x-vector* and *y-vector* entries. This reordering method will be referred to as "mHP$_{\text{RCN}}$", where the small letter "m" is used to indicate the multiple-SpMxV framework.

## 5. Experimental results

### 5.1. Experimental setup

We tested the performance of the proposed methods against three state-of-the-art methods: sBFS [14], sRCM [9, 11, 16] and sHP$_{\text{RN}}$ [18] all of which belong to the single-SpMxV framework. Here, sBFS refers to our adaptation of BFS-based simultaneous data and iteration reordering method of Strout et al. [14] to matrix row and column reordering. Strout et al.'s method depends on implementing breadth-first search on both temporal and spatial locality hypergraphs simultaneously. In our adaptation, we apply BFS on the bipartite graph representation of the matrix, so that the resulting BFS orders on the row and column vertices determine row and column reorderings, respectively. sRCM refers to applying the RCM method, which is widely used for envelope reduction of symmetric matrices, on the bipartite graph representation of the given sparse matrix. Application of the RCM method to bipartite graphs has also been used by Berry et al. [3] to reorder rectangular term-by-document matrices for envelope minimization. sHP$_{\text{RN}}$ refers to the work by Yzelman and Bisseling [18] which utilizes HP using the row-net model for CSR-based SpMxV.

The HP-based top-down reordering methods sHP$_{\text{RN}}$, sHP$_{\text{CN}}$, sHP$_{\text{eRCN}}$ and mHP$_{\text{RCN}}$ are implemented using the state-of-the-art HP tool PaToH [7]. In these implementations, PaToH is used as a 2-way HP tool within the RB paradigm. The hypergraphs representing sparse matrices according to the respective models are recursively bipartitioned into parts until the CSR-storage size of the submatrix (together with the $x$ and $y$ vectors) corresponding to a part drops below the cache size. PaToH is used with default parameters except the use of heavy connectivity clustering (`PATOH_CRS_HCC=9`) in the sHP$_{\text{RN}}$, sHP$_{\text{CN}}$ and sHP$_{\text{eRCN}}$ methods that belong to the single-SpMxV framework, and the use of absorption clustering using nets (`PATOH_CRS_ABSHCC=11`) in the mHP$_{\text{RCN}}$ method that belong to the multiple-SpMxV framework. Since PaToH contains randomized algorithms, the reordering results are reported by averaging the values obtained in 10 different runs, each randomly seeded.

Performance evaluations are carried out in two different settings: cache-miss simulations and actual running times by using OSKI (BeBOP Optimized Sparse Kernel Interface Library) [17]. We provide only the discussion on the OSKI experiments here, whereas we refer the reader to our technical report [1] for the discussion on cache-miss simulations. In OSKI runs, double precision arithmetic is used. OSKI runs are performed on 17 matrices all of which are obtained from the University of Florida Sparse Matrix Collection [15]. CSR-storage sizes of these matrices vary between 13 MB to 94 MB. Properties of these matrices are presented in Table 1. As seen in the table, the test matrices are categorized into three groups as symmetric, square nonsymmetric and rectangular. In each group, the test matrices are listed in the order of increasing number of nonzeros ("nnz"). In the table, "avg" and "max" denote the average and the maximum number of nonzeros per row/column. "cov" denotes the coefficient of variation of number of nonzeros per row/column. The "cov" value of a matrix can be considered as an indication of the level of irregularity in the number of nonzeros per row and column.

### 5.2. OSKI experiments

OSKI experiments are performed by running OSKI version 1.0.1h (compiled with `gcc`) on a machine with 2.66 GHz Intel Q8400 and 4 GB of RAM, where each core pair shares 2 MB 8-way set-associative L2 cache. The Generalized Compressed Sparse Row (GCSR) format available in OSKI is used for all reordering methods. GCSR handles empty rows by augmenting the traditional CSR with an optional list of non-empty row indices thus enabling the multiple-SpMxV framework. For each reordering instance, an SpMxV workload contains 100 calls to `oski_MatMult()` with the same matrix after 3 calls as a warm-up.

Table 2 displays the performance comparison of the existing and proposed methods for the test matrices. In the table, the first column shows OSKI running times without tuning for original matrices. The second column shows the normalized OSKI running times obtained through the full tuning enforced by the `ALWAYS_TUNE_AGGRESSIVELY` parameter for original matrices. The other columns show the normalized running times obtained through the reordering methods. Each normalized value is calculated by dividing the OSKI time of the respective method by untuned OSKI running time for original matrices. As seen in the first two columns of the table, optimizations provided through the OSKI package do not improve the performance of the SpMxV operation performed on the original matrices. This experimental finding can be attributed to the irregularly sparse nature of the test matrices. We should mention that optimizations provided through the OSKI package do not improve the performance of the SpMxV operation performed on the reordered matrices.

Table 1. Properties of test matrices

| Name | number of | | | nnz's in a row | | | nnz's in a column | | |
|---|---|---|---|---|---|---|---|---|---|
| | rows | cols | nonzeros | avg | max | cov | avg | max | cov |
| Symmetric Matrices | | | | | | | | | |
| c-73 | 169,422 | 169,422 | 1,279,274 | 8 | 39,937 | 20.1 | 8 | 39,937 | 20.1 |
| c-73b | 169,422 | 169,422 | 1,279,274 | 8 | 39,937 | 20.1 | 8 | 39,937 | 20.1 |
| rgg_n_2_17_s0 | 131,072 | 131,072 | 1,457,506 | 11 | 96 | 0.3 | 11 | 28 | 0.3 |
| boyd2 | 466,316 | 466,316 | 1,500,397 | 3 | 93,262 | 60.6 | 3 | 93,262 | 60.6 |
| ins2 | 309,412 | 309,412 | 2,751,484 | 9 | 303,879 | 65.3 | 9 | 309,412 | 66.4 |
| rgg_n_2_18_s0 | 262,144 | 262,144 | 3,094,566 | 12 | 62 | 0.3 | 12 | 31 | 0.3 |
| Square Nonsymmetric Matrices | | | | | | | | | |
| Raj1 | 263,743 | 263,743 | 1,302,464 | 5 | 40,468 | 17.9 | 5 | 40,468 | 17.9 |
| rajat21 | 411,676 | 411,676 | 1,893,370 | 5 | 118,689 | 41.0 | 5 | 100,470 | 34.8 |
| rajat24 | 358,172 | 358,172 | 1,948,235 | 5 | 105,296 | 33.1 | 5 | 105,296 | 33.1 |
| ASIC_320k | 321,821 | 321,821 | 2,635,364 | 8 | 203,800 | 61.4 | 8 | 203,800 | 61.4 |
| Stanford_Berkeley | 683,446 | 683,446 | 7,583,376 | 11 | 76,162 | 25.0 | 11 | 249 | 1.5 |
| Rectangular Matrices | | | | | | | | | |
| kneser_10_4_1 | 349,651 | 330,751 | 992,252 | 3 | 51,751 | 31.9 | 3 | 3 | 0.0 |
| neos | 479,119 | 515,905 | 1,526,794 | 3 | 29 | 0.2 | 3 | 16,220 | 15.6 |
| wheel_601 | 902,103 | 723,605 | 2,170,814 | 2 | 442,477 | 193.9 | 3 | 3 | 0.0 |
| LargeRegFile | 2,111,154 | 801,374 | 4,944,201 | 2 | 4 | 0.3 | 6 | 655,876 | 145.9 |
| cont1_l | 1,918,399 | 1,921,596 | 7,031,999 | 4 | 5 | 0.3 | 4 | 1,279,998 | 252.3 |
| degme | 185,501 | 659,415 | 8,127,528 | 44 | 624,079 | 33.1 | 12 | 18 | 0.1 |

As seen in Table 2, on the overall average, the 2D methods sHP$_{eRCN}$ and mHP$_{RCN}$ perform better than the 1D methods sHP$_{RN}$ and sHP$_{CN}$, where mHP$_{RCN}$ (adopting the multiple-SpMxV framework) is the clear winner. Furthermore, for the relative performance comparison of the 1D methods, the proposed sHP$_{CN}$ method performs better than the existing sHP$_{RN}$ method. On the overall average, sHP$_{CN}$, sHP$_{eRCN}$ and mHP$_{RCN}$ achieve significant speedup by reducing the SpMxV times by 11%, 14% and 18%, respectively, compared to the unordered matrices; thus confirming the success of the proposed reordering methods.

Table 3 is introduced to evaluate the preprocessing overhead of the reordering methods. For each test matrix, the reordering times of all methods are normalized with respect to the OSKI time of the SpMxV operation using the unordered matrix and geometric averages of these normalized values are displayed in the "overhead" column of the table. In the table, the "amortization" column denotes the average number of SpMxV operations required to amortize the reordering overhead. Each "amortization value" is obtained by dividing the average normalized reordering overhead by the overall average OSKI time improvement taken from Table 2. Overhead and amortization values are not given for the sRCM method since sRCM does not improve the OSKI running time at all.

As seen in Table 3, top-down HP-based methods are significantly slower than the bottom-up sBFS method. The running times of two 1D methods sHP$_{RN}$ and sHP$_{CN}$ are comparable as expected. As also seen in the table, the 2D methods are considerably slower than the 1D methods as expected. In the column-net hypergraph model used in 1D method sHP$_{CN}$, the number of vertices and the number of nets are equal to the number of rows and the number of columns, respectively, and the number of pins is equal to the number of nonzeros. In the hypergraph model used in 2D methods, the number of vertices and the number of nets are equal to the number of nonzeros and the number of rows plus the number of columns, respectively, and the number of pins is equal to two times the number of nonzeros. That is, the hypergraphs used in 2D methods are considerably larger than the hypergraphs used in 1D methods. So partitioning the hypergraphs used in 2D methods takes considerably longer time than partitioning the hypergraphs used in 1D methods, and the running time difference becomes higher with increasing matrix density in favour of 1D methods. There exists a considerable difference in the running times of two 2D methods sHP$_{eRCN}$ and mHP$_{RCN}$ in favour of sHP$_{eRCN}$. This is because of the removal of the vertices connected by the cut row-nets in the enhanced row-column-net model used in sHP$_{eRCN}$.

As seen in Table 3, the top-down HP methods amortize for larger number of SpMxV computations compared to the bottom-up sBFS method. For example, the use of sHP$_{CN}$ instead of sBFS amortizes after 276% more SpMxV computations on the overall average. As also seen in the table, 2D methods amortize for larger number of SpMxV computations compared to the 1D methods. For example, the use of mHP$_{RCN}$ instead of sHP$_{CN}$ amortizes after 178% more SpMxV computations.

## 6. Conclusion

Single- and multiple-SpMxV frameworks were investigated for exploiting cache locality in SpMxV computations that involve irregularly sparse matrices. For the single-SpMxV framework, two cache-size-aware top-down row/column-reordering methods based on 1D and 2D sparse matrix partitioning were proposed by utilizing the column-net and enhancing the row-column-net hypergraph models of sparse matrices. The multiple-SpMxV framework requires splitting a given matrix into a sum of multiple nonzero-disjoint matrices so that the SpMxV operation is computed as a sequence of multiple input- and output-dependent SpMxV operations. For this framework, a cache-size aware top-down matrix splitting method based on 2D matrix partitioning was proposed by utilizing the row-column-net hypergraph model of sparse matrices. The proposed hypergraph-partitioning (HP) based methods in the single-SpMxV framework primarily aim at exploiting temporal locality in accessing input-vector entries and the proposed HP-based method in the multiple-SpMxV framework primarily aims at exploiting temporal locality in accessing both input- and output-vector entries.

Table 2. OSKI running times for the test matrices (cache size = part-weight threshold = 2 MB)

| | Actual | Normalized w.r.t. Actual Times on Original Order | | | | | | |
| --- | --- | --- | --- | --- | --- | --- | --- | --- |
| | Original | | Existing Methods | | | Proposed Methods | | |
| | Order | | | | | Single SpMxV | | Mult. SpMxVs |
| | not tuned (ms) | OSKI tuned | sBFS [14] | sRCM [11] Modified | sHP$_{RN}$ [18] (1D Part.) | sHP$_{CN}$ (1D Part.) | sHP$_{eRCN}$ (2D Part.) | mHP$_{RCN}$ (2D Part.) |
| Symmetric Matrices | | | | | | | | |
| c-73 | 0.454 | 1.00 | 1.02 | 1.06 | 0.93 | 0.92 | 0.92 | **0.90** |
| c-73b | 0.456 | 1.00 | 1.01 | 1.07 | 0.93 | 0.92 | 0.91 | **0.89** |
| rgg_n_2_17_s0 | 0.503 | 0.95 | 0.92 | 1.07 | 0.89 | 0.82 | **0.76** | 0.91 |
| boyd2 | 0.726 | 1.19 | 1.00 | 1.14 | 0.95 | 0.92 | 0.89 | **0.85** |
| ins2 | 1.207 | 1.00 | 0.96 | 2.32 | 0.97 | 1.06 | 0.97 | **0.67** |
| rgg_n_2_18_s0 | 1.051 | 0.96 | 0.90 | 1.07 | 0.99 | 0.99 | **0.75** | 0.81 |
| Square Nonsymmetric Matrices | | | | | | | | |
| Raj1 | 0.629 | 1.04 | 0.88 | 0.96 | 0.86 | **0.82** | 0.83 | 0.84 |
| rajat21 | 0.953 | 1.07 | 1.01 | 1.16 | 1.00 | 0.95 | 0.96 | **0.90** |
| rajat24 | 0.963 | 1.02 | 1.04 | 1.16 | 0.99 | 0.94 | 0.96 | **0.91** |
| ASIC_320k | 1.436 | 0.99 | 1.09 | 1.44 | 0.97 | 0.92 | 0.73 | **0.64** |
| Stanford_Berkeley | 2.325 | 1.04 | 1.01 | 1.05 | 1.10 | 1.01 | **0.89** | 0.98 |
| Rectangular Matrices | | | | | | | | |
| kneser_10_4_1 | 0.694 | 1.02 | 0.70 | 0.87 | 0.81 | **0.67** | 0.89 | 0.68 |
| neos | 0.697 | 1.26 | 1.14 | 1.19 | 1.00 | **0.95** | **0.95** | 0.96 |
| wheel_601 | 1.377 | 1.27 | 0.82 | 0.75 | 0.69 | 0.69 | 0.66 | **0.52** |
| LargeRegFile | 2.643 | 1.55 | 1.19 | 1.30 | 1.04 | **0.95** | **0.95** | 0.96 |
| cont1_l | 2.939 | 1.14 | 1.04 | 1.19 | 1.05 | **0.93** | **0.93** | 0.95 |
| degme | 2.770 | 1.04 | 0.77 | 1.26 | 0.87 | 0.77 | 0.78 | **0.74** |
| Geometric Means | | | | | | | | |
| Symmetric | - | 1.01 | 0.97 | 1.23 | 0.94 | 0.94 | 0.86 | 0.84 |
| Nonsymmetric | - | 1.03 | 1.00 | 1.14 | 0.98 | 0.93 | 0.87 | 0.84 |
| Rectangular | - | 1.20 | 0.93 | 1.07 | 0.90 | 0.82 | 0.85 | 0.78 |
| Overall | - | 1.08 | 0.96 | 1.15 | 0.94 | 0.89 | 0.86 | 0.82 |

Table 3. Average normalized reordering overhead and average number of SpMxV operations required to amortize the reordering overhead

| | Existing Methods | | Proposed Methods | | | | | | |
| --- | --- | --- | --- | --- | --- | --- | --- | --- | --- |
| | | | Single SpMxV | | | | | | Multiple SpMxVs |
| | sBFS [14] | | sHP$_{RN}$ [18] (1D Part.) | | sHP$_{CN}$ (1D Part.) | | sHP$_{eRCN}$ (2D Part.) | | mHP$_{RCN}$ (2D Partitioning) |
| | Overhead | Amortization | Overhead | Amortization | Overhead | Amortization | Overhead | Amortization | Overhead | Amortization |
| Symmetric | 17 | 465 | 194 | 3135 | 190 | 1716 | 514 | 3732 | 920 | 5097 |
| Nonsymmetric | 26 | 700 | 314 | 5078 | 304 | 2740 | 664 | 4822 | 1198 | 6640 |
| Rectangular | 23 | 621 | 383 | 6197 | 254 | 2292 | 620 | 4503 | 1240 | 6870 |
| Overall | 22 | 587 | 286 | 4620 | 245 | 2209 | 596 | 4327 | 1110 | 6149 |

The performance of the proposed models and methods were evaluated on a wide range of sparse matrices by performing actual runs. Experimental results, which were obtained by performing actual runs through using OSKI, showed that the proposed methods and models outperform the state-of-the-art schemes and also these results conformed to our expectation that temporal locality is more important than spatial locality (for practical line sizes) in SpMxV operations that involve irregularly sparse matrices. The two proposed methods that are based on 2D matrix partitioning were found to perform better than the proposed method based on 1D partitioning at the expense of higher reordering overhead, where the 2D method within the multiple-SpMxV framework was the clear winner.

### Acknowledgements

This work was financially supported by the PRACE project funded in part by the EUs 7th Framework Programme (FP7/2007-2013) under grant agreement no. RI-211528 and FP7-261557. The work is achieved using the PRACE Research Infrastructure resources. We acknowledge that the results in this paper have been achieved using the PRACE Research Infrastructure resource the Gauss Center for Supercomputing IBM Blue Gene/P system JUGENE (Jueliche Blue Gene/P) at Forschungzentrum Juelich (FZJ), Germany.